\documentclass[a4paper,12pt]{article}
\usepackage{amsfonts,amsmath,epsfig,graphicx}
\usepackage{amssymb,exscale,cmmib57}

\newtheorem{thm}{Theorem}
\newtheorem{prop}{Proposition}

\newcommand{\Sy}{\mathrm{S^+(p,n)}}

\newcommand{\RR}{{\mathbb R}}

\newcommand{\SSS}{{\mathbb  S}}

\newcommand{\abs}[1]{\lvert#1\rvert}
\newcommand{\norm}[1]{\lVert#1\rVert}

\newcommand{\dotex}{{\frac{d}{dt}}}

\newcommand{\tr}[1]{\text{Tr}\left(#1\right)}

\usepackage{xcolor}

\title{Riemannian {Metric} and Geometric Mean for Positive Semidefinite Matrices 
of Fixed Rank}

\author{Silv\`{e}re Bonnabel\thanks{Department of Electrical Engineering and Computer 
Science, University of Li\`{e}ge, B-4000 Li\`{e}ge, Belgium (bonnabel@montefiore.ulg.ac.be, 
r.sepulchre@ulg.ac.be).}
      \and Rodolphe Sepulchre\footnotemark[2]}
\date{}
\begin{document}
\maketitle

\begin{abstract}
This paper introduces a new  {metric} and mean on the set of positive 
semidefinite matrices of fixed-rank.  The proposed {metric} is derived from a 
well-chosen  Riemannian quotient geometry that generalizes the reductive 
geometry of the positive cone and the associated natural metric. The resulting 
Riemannian space has strong geometrical properties: it is geodesically complete, 
and the metric is invariant with respect to all transformations that preserve 
angles (orthogonal transformations, scalings, and pseudoinversion).   {A meaningful approximation of the associated Riemannian distance is proposed, that} can be efficiently numerically computed via a simple 
algorithm based on SVD.  The induced mean  preserves the rank, possesses the 
most desirable characteristics of a geometric mean, and is easy to compute.
\end{abstract}




\pagestyle{myheadings}
\thispagestyle{plain}
\markboth{SILV\`{E}RE BONNABEL AND RODOLPHE SEPULCHRE}{METRICS FOR POSITIVE SEMIDEFINITE MATRICES}

\section{Introduction}

Positive definite matrices have become fundamental computational objects in 
many areas of  engineering and applied mathematics. They appear as covariance 
matrices in statistics, as elements of the search space in convex and 
semidefinite programming, as kernels in machine learning, and as diffusion 
tensors in medical imaging, to cite a few. Computing with positive definite 
matrices involves approximations, interpolation, filtering, and estimation, 
leading to a realm of metric-based algorithms. In the recent years, it has been 
increasingly recognized that the Euclidean distance does not suit best the set 
of positive definite matrices---the positive symmetric cone $\mathrm{P_n}$---and 
that working with the proper geometry does matter in computational problems.

The (non-Euclidean) natural metric of the positive cone \cite{faraut} proceeds 
from the rich quotient geometry of this set and its structure of reductive 
homogeneous space (see Appendix~\ref{app2} for more details). The 
resulting ``natural'' metric is invariant  to the action (by congruence) of 
the general linear group, a feature at the core of many desirable properties. 
The associated Riemannian mean can be  called geometric rather than arithmetic. 
Recent contributions that have advocated the use of this metric in applications 
include \cite{pennec-06,fletcher,moakher06} for tensor computing in medical 
imaging and \cite{barbaresco} for radar processing.  More theoretical results 
can be found in \cite{moakher05,ando,arsigny,petz}.  Most notably, the natural 
metric corresponds to the Fisher information metric for the multivariate normal 
distribution \cite{burbea-rao,skovgaard}. The comprehensive paper 
\cite{smith-2005} uses this result to derive an intrinsic Cram\'{e}r--Rao bound 
for the Gaussian covariance matrix estimation problem.  Finally, the natural 
metric coincides with the metric defined by the natural self-concordant 
logarithmically homogeneous barrier\break ($-\log$ det $A$) on the symmetric 
cone (which is a convex set) in optimization \cite{boyd-book,nesterov}. In 
particular, the interest in short-step methods relies on the property that 
$\mathrm{P_n}$ with its natural metric is geodesically complete,   {i.e., every maximal geodesic is defined for all $t\in\RR$, and thus the boundary can not be reached in finite time}.

Because matrix algorithms tend to be applied to computational problems of
ever-increasing size, they need to be adapted to remain tractable. Typical matrix
computations (like SVD, EVD, QR factorization, etc.) require $O(n^3)$ operations for a  positive
definite matrix of size $n$, which limits their use in large-scale problems. A sensible 
remedy is to work with low-rank approximations instead. A rank $p$ approximation
of a positive definite matrix can be factored as $A=Z Z^T$, where the matrix $Z \in \RR^{n \times p}$
is of much reduced size if $p \ll n$,
leading to a reduction of the numerical cost of typical matrix operations from $O(n^3)$ to  $O(np^2)$. If $p$ is kept
to a moderate value, the complexity of the resulting algorithms  grows only
linearly with the size of the problem.

The natural metric developed for positive definite matrices is only valid for
full-rank matrices. The goal of this paper is to extend the natural metric of the positive cone
to the set $\mathrm{S^+(p,n)}$ of symmetric positive semidefinite matrices of fixed-rank $p < n$.
The set $\mathrm{S^+(p,n)}$
admits a quotient geometry that generalizes
the quotient geometry of the cone in a way that preserves most of the desirable properties of the cone geometry.
Motivated by the natural metric of the positive cone, the proposed geometry differs from the 
quotient geometry recently proposed \cite{AbsIshLatHuf2008.013}.  The
resulting  ``natural'' metric---which, to the best of the authors' knowledge, has not appeared in
the literature previously---preserves not all but a remarkably large number of invariance properties of
the natural metric in the cone. More precisely, it is invariant under all transformations that
preserve angles, that is, rotations, scalings, and pseudo-inversion. Endowed with this
metric  $\Sy$ is geodesically complete. {Moreover, we propose a meaningful and numerically efficient approximation of  the Riemannian  distance. The induced notion of mean is } 
is shown to present  all the  desirable properties of a geometric mean.

The proposed natural metric on $\mathrm{S^+(p,n)}$ is viewed as an important 
step to generalize several existing algorithms for positive definite matrices 
to the semidefinite case. This applies not only to the computational problems 
mentioned previously, but also to a growing realm of matrix nearness problems 
based on the use of Bregman divergences. We mention in particular the recent 
paper \cite{dhillon-07} that leaves as an open question the characterization 
of different types of projections to compute distances onto important sets of 
matrices, such as the positive semidefinite cone.

The rest of the paper is organized as follows. In section \ref{PD:sec} we summarize
the existing work on the positive definite cone $\mathrm{P_n}$, concentrating on
the features most relevant for the paper. To gain insight on the main issues
faced when extending the geometry of $\mathrm{P_n}$ to $\mathrm{S^+(p,n)}$, section \ref{2D:sec} focuses
on  the simplest case of $2\times 2$ matrices of rank 1.
A metric extending the natural metric on the cone is derived using
polar decomposition. We  provide
geometrical and physical justifications for the associated invariance, distance,
and mean properties.

Building upon the polar decomposition of vectors in the plane, we develop in section \ref{dist:def:sec} 
a quotient  geometry for   $\Sy$. The resulting
natural metric decomposes as the sum of the natural metric on the cone and the
standard metric of the Grassman manifold. We prove that $\Sy$ endowed with
this metric is a Riemannian manifold.

Section~\ref{dist:sec} is devoted to { the Riemannian distance associated
to the metric. Because the explicit calculation of geodesics is out of reach, we construct special curves that approximate the geodesics. The length of those curves can be calculated by a SVD based algorithm of numerical complexity $O(np^2)$. It provides a meaningful  notion of closeness in $\Sy$ which inherits  the invariance properties of the Riemannian distance. The singularities
of this measure are also characterized.  As an aside, the manifold $\Sy$ is proved to be geodesically complete.  }

Section~\ref{mean:sec} provides a {new definition of geometric} mean between two
matrices of $\Sy$ {based on the measure of closeness. We prove that the mean preserves the rank, we 
argue it deserves the appellation ``geometric'', and we prove it is a generalization to rank-deficient positive symmetric matrices of the mean associated with the natural metric on the positive cone.} From a computational viewpoint, a main advantage of the 
metric is that, after a suitable SVD, the computations of the {closeness} and mean 
decouple into two separate problems: computing a Riemannian distance and mean 1. 
in the (lower dimension) cone $\mathrm{P_p}$ and 2. in the Grassman manifold of 
subspaces of dimension $p$ in $\RR^n$.

Concluding remarks and perspectives are discussed in section \ref{s7}.

\subsection{Notation}
\begin{itemize}
  \item $\mathrm{P_n}$ is the set of symmetric positive definite $n\times n$ matrices.
  \item $\mathrm{S^+(p,n)}$ is the set of symmetric positive semidefinite $n\times n$ matrices of rank $p\leq n$. 
  We will only use this notation in the case $p<n$.
  \item $\mathrm{Gl(n)}$ is the general linear group, that is, the set of invertible $n \times n $ matrices.
   \item $\RR_*^{n\times p}$ is the set of full rank $n\times p$ matrices.
   \item $\mathrm{V_{n,p}}=\mathrm{O(n)}/\mathrm{O(n-p)}$ is the  Stiefel manifold; i.e., the set of $n\times p$
matrices with orthonormal columns: $U^TU=I_p$.
\item $\mathrm{Gr(p,n)}$ is  the Grassman manifold, that is,  the set of $p$-dimensional subspaces of $\RR^n$. 
It can be represented by the equivalence classes $\mathrm{V_{n,p}}/\mathrm{O(p)}$.
  \item $\mathrm{Sym(n)}$ is the vector space of symmetric $n\times n$ matrices.
  \item diag$(\lambda_1,\ldots,\lambda_n$) is the $n\times n$ matrix with the $\lambda_i$'s on its diagonal. 
  $I={\rm diag}(1,\ldots,1)$ is the identity matrix.
  \item range($A$) is the subspace of $\RR^n$ spanned by the columns of $A\in \RR^{n\times n}$.
  \item $T_X\mathcal M$ is the tangent space to the manifold $\mathcal M$ at $X$.
\end{itemize}

\section{Riemannian  distances and geometric means on the symmetric cone}\label{PD:sec}
The geometry of the $n$-dimensional symmetric cone $\mathrm{P_n}$  has been well-studied in the literature.
This section reviews some of its relevant features in view of the main developments of the present paper.
Given a matrix $A\in \mathrm{P_n}$, a starting point is the matrix factorization
\begin{align}\label{matrix:fac}
A=ZZ^T=(UR)(UR)^T=UR^2U^T=R'^2,
\end{align}
where $A\in \mathrm{P_n}$, $Z\in \mathrm{Gl(n)}$, $R, R' \in \mathrm{P_n}$, $U\in \mathrm{O(n)}$.
The right and left polar decompositions $Z=UR=R' U $ are unique \cite{golub-book}, but the factorization $A=Z Z^T$
is unaffected by orthogonal transformations $Z\mapsto Z O$ with $ O \in \mathrm{O(n)}$. The matrix equalities \eqref{matrix:fac} underline the quotient geometry of the cone $\mathrm{P_n}$.
\begin{align}\label{quot:cone}
P_n\cong \mathrm{Gl(n)}/\mathrm{O(n)}\cong (\mathrm{O(n)}\times \mathrm{P_n})/\mathrm{O(n)}.
\end{align}
The characterization \eqref{quot:cone} encodes the rich geometry of $\mathrm{P_n}$ as a reductive homogeneous space,
as discussed in \cite{jost-book,smith-2005}.  The most relevant consequence of this feature to the present context is the existence
of a $\mathrm{Gl(n)}$-invariant metric on the manifold $\mathrm{Gl(n)}/\mathrm{O(n)}$. This metric is called the natural metric
on the symmetric cone \cite{faraut}. Up to a scaling factor it is also known as the affine-invariant metric \cite{pennec-06},
the Siegel  metric in symplectic geometry, and it coincides with the metric given by the Fisher information matrix for Gaussian covariance matrix estimation \cite{smith-2005}.
We briefly summarize how it is derived and its main properties in the present context.

In view of \eqref{matrix:fac}, $\mathrm{Gl(n)}$ has a transitive action on $\mathrm{P_n}$  via
congruence
\begin{align}\label{group:action:def}
A\mapsto LAL^T,
\end{align}
where $L\in \mathrm{Gl(n)}$ and any matrix $A$ is brought back to the identity matrix
choosing $L=A^{-1/2}$ with $A^{1/2}$ defined by the polar factor $R'$ in \eqref{matrix:fac}.
Likewise, any tangent vector $X\in \mathrm{Sym(n)}$ at identity $I\in \mathrm{P_n}$ can be transported to a tangent vector $A^{1/2}XA^{1/2}\in T_A P_n$. At identity $I$, the
$\mathrm{Gl(n)}$-invariant metric is defined as the usual scalar product
$g^{P_n}_I(X_1,X_2)=\tr{X_1X_2^T}=\tr{X_1X_2}$. The invariance of the metrics then implies $g^{P_n}_A(A^{1/2}X_1A^{1/2},A^{1/2}X_2A^{1/2})=\tr{X_1X_2}$, which can only be satisfied with the definition
\begin{align}\label{siegel:metr:def}g^{P_n}_A(D_1,D_2)=\tr{D_1A^{-1}D_2A^{-1}}\end{align}at any arbitrary $A\in \mathrm{P_n}$.
The invariance of the metric \eqref{siegel:metr:def} has direct implications on the expression of the geodesics and the accompanying Riemannian distance.
The exponential map at $I$ is the usual matrix exponential
\[
\exp^{P_n}_IX=\exp X=\sum_{k=0}^\infty({1}/{k!})X^k.
\]
The Frobenius  norm $\norm{X}_F$ is the geodesic length
$d(\exp X,I)=\norm{X}_F$, hence the formula at identity $$d(A,I)=\norm{\log A}_F.$$ Invariance of 
the metric again extends the characterization of geodesics at arbitrary
$A\in \mathrm{P_n}$ \cite{moakher05,smith-2005}:$$\exp^{P_n}_A(tX)=A^{1/2}\exp(tA^{-1/2}XA^{-1/2})A^{1/2},\quad
t>0,$$and the corresponding geodesic distance\begin{align}\label{sieg:dist:frob}d_{P_n}(A,B)=
d(A^{-1/2}BA^{-1/2},I)
=\norm{\log(A^{-1/2}BA^{-1/2})}_{F}=\sqrt{\sum_k \log^2(\lambda_k)},\end{align}where $\lambda_k$ are
the generalized eigenvalues of the pencil $A-\lambda B$, i.e., the roots of
$\det(AB^{-1}-\lambda I)$. Note that the distance is invariant with respect to matrix inversion $(A,B)\mapsto(A^{-1},B^{-1})$ because $\log^2(\lambda_k)$
is invariant to inversion $\lambda_k\mapsto \lambda_k^{-1}$.

The geodesic characterization provides a closed-form expression of the Riemannian (Karcher) 
mean of two matrices $A,B\in \mathrm{P_n}$. The geodesic $A(t)$ linking $A$ and $B$ is 
$$A(t)=\exp^{P_n}_A(tX)=A^{1/2}\exp(t\log(A^{-1/2}BA^{-1/2}))A^{1/2},$$ where 
$A^{-1/2}XA^{-1/2}=\log(A^{-1/2}BA^{-1/2})\in \mathrm{Sym(n)}$. The midpoint is obtained 
for $t=1/2$:\begin{align}\label{geom:classical:mean}A\circ B=A^{1/2}(A^{-1/2}BA^{-1/2})^{1/2}A^{1/2}.\end{align}

For a comprehensive treatment of the geometric versus arithmetic means of positive 
definite matrices, the reader is referred to \cite{ando,moakher05}. From a geometric 
viewpoint, there are numerous reasons to prefer the natural metric (and associated 
distance and mean) on the symmetric cone to the flat metric associated to the distance 
$\norm{A-B}_F$. Indeed the symmetric cone is $not$ a vector space. The flat metric can 
still be used because the set is convex. However,  many applications treating the space 
of covariance matrices as a vector space can yield to degraded algorithm performances 
\cite{smith-2005}. A further undesirable feature of the flat metric is that 
$\mathrm{P_n}$ is not a geodesically complete space (see section \ref{dist:sec}) since 
the geodesic $A+t(B-A)$ is not a positive matrix for all $t$. In contrast, it becomes 
geodesically complete with the natural metric. A practical consequence is that the 
natural metric is well-suited to short-step interior point methods in $\mathrm{P_n}$. 
Remarkably, it coincides with the metric defined by the natural self-concordant 
logarithmically homogeneous barrier ($-$log det $A$) on the symmetric cone (which is a 
convex set) in optimization \cite{boyd-book,nesterov}. Furthermore, the invariance to 
the group action \eqref{group:action:def} (implying in particular invariance with 
respect to inversion) is desirable in numerous applications (see, e.g., 
\cite{arsigny,pennec-06,smith-2005,barbaresco}). In particular, if $A$ is a covariance 
matrix $\mathbb{E}(zz^T)$, the action corresponds to a change of basis $z\mapsto Lz$.  
As a consequence the natural metric is well-suited to intrinsic estimation algorithms 
for covariance matrix estimation; see \cite{smith-2005}. In 
\cite{pennec-06,fletcher,arsigny}, the authors argue that invariance with respect to 
inversion is adapted to the physics of diffusion tensors related to medical imaging. 
Finally, the associated (geometric) mean has numerous desirable properties developed 
in section \ref{mean:sec}. It is useful to keep in mind that the geometric mean 
coincides with the classical geometric mean $\sqrt{ab}$ when $a,b\in \mathrm{P_1}$.

\section{Extending the metric: A geometric insight in the plane}\label{2D:sec}
In order to inherit some of the nice properties of the natural
distance on $\mathrm{P_n}$ when the  matrices are not full-rank, we
seek a distance which is invariant to the $\mathrm{Gl(n)}$ group action. For tutorial purposes we start by
considering $2\times 2$ matrices of rank $1$: $\mathrm{S^+(1,2)}$. For any
$L\in \mathrm{Gl(2)}$ and $A,B\in \mathrm{S^+(1,2)}$ one wants
$$
d(A,B)=d(LAL^T,LBL^T).
$$
{This objective is, however, too ambitious}.
Taking
$$A=\begin{pmatrix}1&0\\0&0\end{pmatrix},\quad  B=\begin{pmatrix}1&1\\1&1\end{pmatrix},\quad \text{and} \quad L=\begin{pmatrix}1&0\\0&\epsilon\end{pmatrix},$$
$A$ is unchanged by the transformation whereas
$LBL^T=\left(\begin{smallmatrix}1&\epsilon\\\epsilon&\epsilon^2\end{smallmatrix}\right)$.
Thus if the distance is continuous in the matrix elements (for example it is associated 
to a metric), we have in the limit $ d(A,B)=d(A,A)=0 $ as $\epsilon\rightarrow 0$. This 
proves that $d$ is not a distance since $A\neq B$.

Mimicking the matrix decomposition \eqref{matrix:fac}, write $A\in
\mathrm{S^+(1,2)}$: $$A=xx^T=ur^2u^T,$$ where $u=(\cos (\theta),
\sin (\theta)) \in\SSS^1\subset\RR^2$ is a unit vector and
$(r,\theta)$ is the polar representation of $x$. Without loss of
generality  $\theta$ is equated to $\theta+j\pi$, $j\in \mathbb
Z$, since $x$ and $-x$ correspond to the same $A$. Thus
$\mathrm{S^+(1,2)}$ can be equated to the space $\RR^{+}_*\times
\RR\mathbb P^1$, where $\RR\mathbb P^1$ is the real projective
space of dimension 1 (lines of $\RR^2$). The group action
\eqref{group:action:def} corresponds to a change of basis for the
vector $x\mapsto Lx$. We have proved in the last paragraph that
there is no distance which is invariant under every change of basis
for $x$. Nevertheless one can define a distance which is invariant
by scaling and orthogonal transformation, i.e., $G=\{\mu R:
(\mu,R)\in \RR\times O(2)\}$.  A sensible metric is $ds^2=
d\theta^2+{k}~(dr/r)^2= d\theta^2+{k}~d(\log r)^2$ {with $k>0$}. Let us find the
distance associated to this metric. Let $A=x_1x_1^T$ and
$B=x_2x_2^T$. Let $(r_i,\theta_i)$ denote their
polar 
coordinates for $i=1,2$. It is always
possible to have $\abs{\theta_2-\theta_1}\leq \pi/2$ possibly
replacing $\theta_1$ by $\theta_1\pm\pi$. The distance  is then
\begin{align}\label{metr:vec:def}
d_{\mathrm{S}^+(1,2)}^2(A,B)=\abs{\theta_2-\theta_1}^2+k~\abs{\log({r_1}/{r_2})}^2,\quad {k>0}.
\end{align}
where the first term penalizes the distance between the
subspaces range($A$)= span($x_1$) and range($B$)=span($x_2$). It ensures the invariance by rotation/orthogonal
transformation (which only affects $\theta$) while the second term ensures the invariance  by scaling (which only affects $r$). Note that when
range($A$)=range($B$) the induced distance on the 1-dimensional
subspace corresponds to the natural distance \eqref{sieg:dist:frob}
between matrices of $P_1$. The metric above rewrites
in the general form
\begin{align}\label{metr:case2:def}ds^2= \tr{du^Tdu}+{k}~d(\log
r)^2,\quad {k>0}.\end{align}

\paragraph{A symmetry-based justification}
The choice of the metric \eqref{metr:case2:def} can be derived
from necessary conditions imposed by  the desired symmetries (invariances).
Let $A=x_1x_1^T$, $B=x_2x_2^T$ with $x_1,x_2\in\RR^2$. We seek
an invariant distance between $A$ and $B$. This is equivalent to seeking an invariant 
distance between $x_1$ and $x_2$; i.e., a scalar $G$-invariant function on $\RR^2\times \RR^2$,
satisfying the conditions required  to be a distance. If $H$ is any subgroup of $G$, 
the distance must be in particular $H$-invariant. Let $H=\RR^+_*\times SO(2)\subset G$.
The identification $\RR^2_*\cong H$ can be obtained via polar coordinates 
$(r,\theta)\in \RR^+_*\times S^1\cong H$. A standard result  (see, e.g., \cite{olver-book95,arxiv-07})
is that  every $H$-invariant scalar function of ($x_1,x_2)\in H\times H$ (in particular any distance)
is a function of  the scalar invariants $(r_2/r_1,\theta_2-\theta_1)$; i.e., the 
coordinates of $x_1^{-1}*x_2$, where * is the group multiplication of $H$.
Thus any invariant distance $d$ writes  $d(x_1,x_2)=h(r_2/r_1,\theta_2-\theta_1)$, 
and since $d$ is a distance it must be symmetric $d(x_1,x_2)=d(x_2,x_1)$.
Thus necessarily $d(x_1,x_2)=f(\abs{\log(r_2/r_1)}, \abs{\theta_2-\theta_1})$ with 
$f$ positive, monotone in each of its arguments, and only $f(0,0)$ is equal to $0$.
It proves that the distance \eqref{metr:vec:def} is a prototype for every $G$-invariant distance $d(A,B)$.

\paragraph{A physical  justification}The invariance properties of the distance \eqref{metr:vec:def} are meaningful from a physical viewpoint.
Suppose  $Z=\mathbb E(zz^T)$ is the covariance matrix of a stochastic variable $z\in\RR^2$ (for instance the position of an object) with zero mean
such that every realization of $z$ is on a line. Then $Z\in \mathrm{S^+(1,2)}$. The distance \eqref{metr:vec:def} between the two covariances $Z_1, Z_2$
of two independent variables $z_1,z_2$ (dispersing in two different directions) is well-defined. It does not depend on any nontrivial way on the choice of measurement
units, e.g., feet versus meters, as well as on the orientation of the
frame chosen, e.g.,  the first axis is pointing north or south.

Moreover, suppose the measurements are noisy. For instance 
$Z_1=\text{diag(}4,\epsilon^2)$ and $Z_2=\text{diag}(\epsilon^2,1)$, where 
the term  $\epsilon \ll 1$ is the amplitude of the covariance of the noise. 
These two matrices belong to $P_2$, and the geometric mean 
\eqref{geom:classical:mean} is  diag($2\epsilon , \epsilon$). The smallness 
of the noise ruins the mean, which no longer reflects the physical 
interpolation between the two processes. Indeed in the degenerate case 
$\epsilon\rightarrow 0$, the mean becomes the null matrix. In contrast, the 
midpoint in the sense of the distance \eqref{metr:vec:def} between the rank 1 
approximations of the two matrices (i.e., $\epsilon=0$) is
$\left(\begin{smallmatrix}1&1\\1&1\end{smallmatrix}\right)$ with 
eigenvalues $(0,2)$, which is an image of the covariance matrix of the 
midpoint. The distance between the matrices and their mean is independent 
of the choice of units and orientation, and is hardly affected by noise.

\section{A new Riemannian metric on the symmetric semidefinite cone}\label{dist:def:sec}

Let $\mathrm{S^+(p,n)}$ be the set of positive semidefinite matrices of 
fixed-rank $p<n$. One can prove (analogously to the last section) in the 
general case that is impossible to find a distance between matrices of 
$\mathrm{S^+(p,n)}$ which would be invariant to the transformation
\eqref{group:action:def} for $p$ arbitrary small. Nevertheless one can construct 
a metric generalizing \eqref{metr:case2:def} and a distance generalizing 
\eqref{metr:vec:def}. The Grassmanian distance is a multidimensional 
generalization of the angular distance $d\theta$, and the natural metric on 
the symmetric cone \eqref{siegel:metr:def} is a multidimensional 
generalization of $d\log r=r^{-1}dr$. The resulting metric is invariant to 
orthogonal transformations and scalings.   {Locally,} when the ranges of two {infinitely close} matrices coincide, the  induced
metric on the corresponding  subspace {reduces to} the natural metric
on the symmetric cone $\mathrm{P_p}$. In particular, when $p=n$ 
the metric coincides with the natural metric on $\mathrm{P_n}$.

Mimicking the developments of section \ref{PD:sec}, we start from the matrix equalities
$$
A=ZZ^T=(UR)(UR)^T=UR^2U^T,
$$
where $A\in \Sy$, $Z\in \RR_*^{n\times p}$, $R\in \mathrm{P_p}$, $U\in \mathrm{V_{n,p}}$.
Right multiplication of $Z$ by an orthogonal matrix $Z\mapsto ZO,~O\in \mathrm{O(p)}$, does 
not affect the product $A=ZZ^T$. Consider the following $\mathrm{O(p)}$ group action:\begin{align*}
R&\mapsto O^TRO\quad\in \mathrm{P_p},\\
U&\mapsto UO\quad\in \mathrm{V_{n,p}}.
\end{align*}The representation $A=UR^2U^T$ with $(U,R^2)\in \mathrm{V_{n,p}}\times \mathrm{P_p}$ 
is thus univocal up to the equivalence relation $(U,R^2)\equiv (UO,O^TR^2 O)$ 
for any $O\in \mathrm{O(p)}$. Thus the set $\Sy$ admits a quotient manifold representation 
$$\Sy\cong(\mathrm{V_{n,p}}\times \mathrm{P_p})/\mathrm{O(p)}.$$ (A dimension checking yields that
the dimension of $\Sy$ is dim($\mathrm{V_{n,p}}\times \mathrm{P_p}$)-dim($\mathrm{O(p)})=pn-p(p-1)/2$.)
Note that the chosen quotient geometry differs from the one recently considered in 
\cite{AbsIshLatHuf2008.013}, where elements of $\Sy$ are represented
by equivalences classes $ZO,~O\in \mathrm{O(p)}$, leading to the quotient representation 
$\Sy \cong \RR_*^{n\times p} /\mathrm{O(p)}$ .

If $(U,R^2)\in \mathrm{V_{n,p}}\times \mathrm{P_p}$ represents $A\in\Sy$, it is tempting to represent the tangent vectors of $T_A\Sy$ by the infinitesimal variation $(\Delta, D)$, where
\begin{equation}\label{tangent:space}
\begin{aligned}\Delta&=U_\perp B,\quad\quad B\in \RR^{(n-p)\times p},\\D&=RD_0R\end{aligned}
\end{equation} such that $U_\perp\in V_{n,n-p}$ , $U^TU_\perp=0$, and $D_0\in \mathrm{Sym(p)}=T_I\mathrm{P_p}$. The chosen metric of $\Sy$ is merely
the sum of the infinitesimal distance in $\mathrm{Gr(p,n)}$ and  in $\mathrm{P_p}$:
\begin{align}\label{metric:def} g_{(U,R^2)}((\Delta_1,D_1),(\Delta_2,D_2))=\tr{\Delta_1^T\Delta_2}+{{k}}~\tr{R^{-1}D_1R^{-2}D_2R^{-1}},\quad {{k}}>0,
\end{align}generalizing \eqref{siegel:metr:def} in a natural way. The next theorem proves that the construction endows the space $\Sy$ with a Riemannian structure.
\begin{thm}\label{metric:thm}
The space $\Sy\cong (\mathrm{V_{n,p}}\times \mathrm{P_p})/\mathrm{O(p)}$ endowed with the metric \eqref{metric:def} is a Riemannian manifold with horizontal space $$\mathcal
H_{(U,R^2)}=\{(\Delta,D): \Delta=U_\perp B, ~B\in \RR^{(n-p)\times p},
~D=RD_0R,~D_0\in \mathrm{Sym(p)}\}.$$
Furthermore, the metric is invariant with respect to orthogonal transformations, scalings, and pseudoinversion.
\end{thm}
\unskip

\noindent Proof:
In this proof we also recap some results on quotient manifolds. We
follow the machinery of Riemannian quotient manifold (see, e.g.,
\cite{absil-book}), except that we will not require that the quotient map be a Riemannian submersion.
Any representative $(U,R^2)$ lives in the structure space $\mathrm{V_{n,p}}\times
\mathrm{P_p}$. The tangent space to $(U,R^2$) is the direct sum of a vertical
space and a horizontal space: $T_{(U,R^2)}(\mathrm{V_{n,p}}\times \mathrm{P_p})=\mathcal
V_{(U,R^2)}\oplus\mathcal H_{(U,R^2)}$. The equivalence class
$\{(UO,O^TR^2O), ~O\in \mathrm{O(p)}\}$ is called a fiber. The vertical
space is the tangent space to the fiber at $(U,R^2)$:$$\mathcal
V_{(U,R^2)}=\{(U\Omega,R^2\Omega-\Omega R^2): \Omega \in T_I\mathrm{O(p)}\},$$
i.e.,
$\Omega$ is a skew-symmetric matrix. The horizontal space $\mathcal
H_{(U,R^2)}$ is by definition complementary to $\mathcal V_{(U,R^2)}$.
It is customary to represent tangent vectors to the quotient
manifold only as elements of the horizontal space. Indeed it is
useless to consider elements of the vertical space since they are
tangent to the fiber, and all elements of the fiber represent the
same point in the quotient manifold.  Here  the tangent space to the
product manifold writes $T_{(U,R^2)}(\mathrm{V_{n,p}}\times
\mathrm{P_p})=\{(\Xi,D):~~\Xi\in T_U \mathrm{V_{n,p}},~D\in T_{R^2} \mathrm{P_p}\}$. So $\Xi=\Delta
+U\Omega$ with $\Omega$ skew-symmetric and $\Delta=U_\perp B$. But
$(U\Omega,R^2\Omega-\Omega R^2)\in\mathcal V_{(U,M)}$. So the
horizontal space at $T_{(U,R^2)}$ can be chosen to be made only of the
vectors $(\Delta,D$) given by \eqref{tangent:space}: $$\mathcal
H_{(U,R^2)}=\{(\Delta,D): \Delta=U_\perp B, ~B\in \RR^{(n-p)\times p},
~D=RD_0R,~D_0\in sym(p)\}.$$
The manifold $\Sy$ endowed
with the metric \eqref{metric:def} is a Riemannian quotient
manifold: after having chosen  a Riemannian metric on the structure
space, all we must prove is that
the induced metric on the horizontal space does \emph{not} depend on the representative chosen. The Riemannian metric for the product manifold (structure space)
$\mathrm{V_{n,p}}\times \mathrm{P_p}$ can be chosen as the sum of the natural metrics of $\mathrm{V_{n,p}}$ and $\mathrm{P_p}$:
$g_{(U,R^2)}^{SP}(X_1,X_2)=\tr{\Omega_1^T\Omega_2+\Delta_1^T\Delta_2}+{k}~\tr{D_1R^{-2}D_2R^{-2}}$, where $X_i=(U\Omega_i+\Delta_i,D_i)$ for $i=1,2$.
Let $X_1,X_2\in\mathcal H_{(U,R^2)}$. We have $g_{(U,R^2)}^{SP}(X_1,X_2)=\tr{\Delta_1^T\Delta_2}+k~\tr{D_1R^{-2}D_2R^{-2}}$. The first term only depends
on $U_\perp$, which is invariant along the fiber, and the second term is $\tr{D_{0,1}D_{0,2}}$, where $D_i=RD_{0,i}R$, and it does not depend on the representative chosen.
These invariances are due indeed to the invariance of the Grassman metric with respect to the representative in the Stiefel manifold, and  the invariance properties
of the natural metric on $\mathrm{P_p}$.

Concerning the last point of the theorem, an orthogonal transformation affects the representative $(U,R^2)$ by transforming  $U$ in $OU$ with $O\in \mathrm{O(n)}$
and thus $\Delta$ in $O\Delta$; and a scaling by transforming $R^2$ in $\mu^2R^2$ with $\mu \in \RR_*$ and thus $D$ in $\mu^2D$. Thus they affect separately
the first and second term of the metric, which are both invariant to each of these transformations. The invariance to pseudoinversion derives from the invariance
to inversion of the natural metric on $\mathrm{P_p}$, as $(U,R^{-2})$ is a representative of a pseudoinverse of $UR^2U^T$. These properties are further detailed
in the proof of Theorem~\ref{inv:dist:thm}.\qquad

The reader will note that with the proposed choice of metric and horizontal space, $\mathcal V_{(U,R^2)}$ and $\mathcal H_{(U,R^2)}$ are complementary,
but not orthogonal. {This is rather unusual, and it will have implications in the sequel (Proposition 1). However, they tend to be orthogonal as $k\rightarrow 0$}. 

As a concluding remark for this section, we recall that the natural metric on the cone is not the only  $\mathrm{GL(n)}$-invariant metric,
and a whole family of invariant metrics can be derived from the scalar product at identity $\tr{D_1^TD_2}+\beta\tr{D_1}\tr{D_2}$, with $\beta>-\frac{1}{n}$
(see Appendix \ref{app2}). For instance, the metric corresponding to $\beta=-\frac{1}{n+1}$ was proposed in \cite{Lovric} as a Riemannian metric
on the space of Gaussian distributions having a  nonzero mean. This family of metrics could also easily be extended to $\Sy$. It should be clear that the results
of Theorem~\ref{metric:thm} hold for any $\mathrm{GL(n)}$-invariant metric on the cone and not only for the natural one.

\section{Horizontal geodesics in the structure space approximate geodesics in \boldmath$\Sy$} \label{dist:sec} This section provides the construction of special curves of interest connecting any two matrices  $A$ and $B$ in  $\mathrm{S^+(p,n)}$ : horizontal geodesics in the structure space. First of all let us find two representatives of $A$ and $B$ in $\mathrm{V_{n,p}}\times\mathrm{P_{p}}$ connected by a horizontal geodesics. Let $V_A,V_B\in \mathrm{V_{n,p}}$ be two matrices that span range($A$) and range($B$), respectively. The SVD of $V_B^T V_A$ yields $O_A,O_B\in \RR^{p\times p}$  such that \begin{align}\label{SVD}
O_A^TV_A^T V_BO_B=\text{diag}(\sigma_1,\cdots,\sigma_p), \quad 1\geq
\sigma_1 \geq \cdots\geq \sigma_p \geq 0
\end{align}The $\sigma_i=\cos{\theta_i}$ are the cosines of the principal angles $0\leq \theta_1 \leq..\leq\theta_p\leq\pi/2$ between the two subspaces \cite{golub-book}. Choosing the principal vectors  $U_A=(u^A_1,\cdots,u^A_p)=V_AO_A$ and $U_B=(u_1^B,\cdots,u^B_p)=V_BO_B$ yields a simple formula for the Grassman geodesic connecting range($A$) and range($B$) (e.g. \cite{smith-2005}): \begin{align}\label{motion:eq}U(t)=U_A\cos{(\Theta t)}+X\sin(\Theta t)\end{align}where
$\Theta=\text{diag}(\theta_1,\cdots,\theta_p)$ and $X$ is the
normalized projection of $V$ onto the column space of $U_\perp$, i.e., $X=(I-U_AU_A^T)U_BF$ where $F$
is the pseudoinverse of the matrix diag$(\sin(\theta_1),\cdots,\sin(\theta_p))$. The associated geodesic $R^2(t)$ in $\mathrm{P_p}$ must connect $R_A^2:=U_A^TAU_A$ and $R_B^2:=U_B^T BU_B$, that is, 
\begin{align}\label{conegeodesic:eq}
R^2(t)=R_A\exp(t\log R_A^{-1}R_B^2R_A^{-1})R_A.
\end{align}
\begin{thm}\label{inv:dist:thm}  The singular value decomposition \eqref{SVD} and the geodesic curves (\ref{motion:eq}) and (\ref{conegeodesic:eq}) define a curve in $\Sy$
\begin{align}\label{curve:def}
\gamma_{A \rightarrow B}(t)=U(t)R^2(t)U^T(t)\end{align}
with the following properties:
\begin{itemize}
\item $\gamma_{A \rightarrow B}(\cdot)$ connects $A$ and $B$ in $\mathrm{S^+(p,n)}$, that is,  $\gamma_{A \rightarrow B}(0)=A$,  $\gamma_{A \rightarrow B}(1)=B$, and $\gamma_{A \rightarrow B}(t) \in \mathrm{S^+(p,n)} ~\forall t \in [0,1]$.
\item The curve $(U(t),R^2(t))$ is a horizontal lift of $\gamma_{A \rightarrow B}(t)$ and it is a geodesic in the structure space $V_{n,p} \times P_p$. 
\item The (squared) total length of $\gamma_{A \rightarrow B}(t)$ in the Riemannian manifold $(\mathrm{S^+(p,n), g)}$ is given by 
\begin{align}\label{dist:def}l^2(\gamma_{A \rightarrow B})&=\norm{\Theta}_F^2+k~\norm{\log R_A^{-1}R_B^2 R_A^{-1}}_F^2 . \end{align} 
It is invariant with respect to pseudoinversion and to the group action by congruence of orthogonal transformations and scalings.
\end{itemize}
Furthermore, the curve $\gamma_{A \rightarrow B}(\cdot)$ is uniquely defined provided that the $(p-1)$th principal angle satisfies $\theta_{p-1}\neq \pi/2$.
\end{thm}

The proof of the theorem is given in Appendix \ref{app2}. Viewing matrices of $\Sy$ as flat ellipsoids in $\RR^n$, the length  $l(\gamma_{A \rightarrow B})$ consists of two independent contributions: a distance between the subspaces in which the ellipsoids are contained (Grassman distance), and a distance between the ellipsoids within a common subspace (natural distance on the cone). In this sense, \eqref{dist:def} provides an exact generalization of the Riemannian distance \eqref{metr:vec:def}  of Section 3. However, it is {\it not} the Riemannian distance of $(\mathrm{S^+(p,n), g)}$ because $\gamma_{A \rightarrow B}(\cdot)$ is not necessarily a geodesic curve (even tough it is the base curve of a horizontal geodesic in the structure space). 
\begin{prop}
The curve $\gamma_{A \rightarrow B}(\cdot)$ is not necessarily a geodesic of $(\mathrm{S^+(p,n), g)}$. Its length provides a meaningful measure of {\it closeness} between $A$ and $B$ which is not a distance because it does not satisfy the triangle inequality.
\end{prop}

A proof of this proposition is given by the following example:   Let $A=\text{diag}(2,1,0)$ and  $B=\text{diag}(1,2,0)$.  Because $A$ and $B$ have the same range, the length  $l(\gamma_{A \rightarrow B})$ reduces to a distance in the cone $P_2$ (no Grassman contribution). But $A$ can also be connected to $B$ via $C=\text{diag}(0,1,2)$ and $D=\text{diag}(1,0,2)$. The curves  $\gamma_{A \rightarrow C}(\cdot)$, $\gamma_{C \rightarrow D}(\cdot)$, and $\gamma_{D \rightarrow B}(\cdot)$ each involve a pure subspace rotation of $\pi/2$, which means that their total length in  $(\mathrm{S^+(p,n), g)}$ is $3\pi/2$ (no contribution in the cone). The situation $l(\gamma_{A \rightarrow B})> 3\pi/2$ is clearly possible if $k$ is large enough, showing that $\gamma_{A \rightarrow B}(\cdot)$ is not necessarily a geodesic curve. Even if the triangular inequality is not satisfied, the proposed measure of closeness satisfies the two other properties of a distance: it is symmetric, and vanishes only when the two matrices coincide.

The construction of the curve  $\gamma_{A \rightarrow B}(\cdot)$   has a geometrical meaning in $\Sy$: If a gyroscope is attached to the moving ellipsoid represented by $U(t)R_A^2U(t)$, where $U(t)$ is given by \eqref{motion:eq}, it indicates no rotation around an axis perpendicular to the ellipsoid during the motion $0\leq t\leq 1$. Thus the ellipsoid $U(1)R_A^2U(1)=U_BR_A^2U_B$ is the ellipsoid $A=U(0)R_A^2U(0)=U_AR_A^2U_A$ brought in range($B$) by a rotation of minimal energy. This also justifies  comparing directly $R_A^2$ to $R_B^2$ in $\mathrm{P_p}$. If $O\neq I$ is an orthogonal matrix, it would violate the rotational invariance to compare $R_A^2$ to $O^TR_B^2O$, even though $(U_BO,O^TR_BO)$ is a valid representation of $B$ and the distance in Grassman between range($U_A$) and range($U_B$) is unaffected by the transformation $U_B\mapsto U_BO$.

Even tough it is not a distance, the closeness  $l(\gamma_{A \rightarrow B})$  is a meaningful generalization of the Riemannian distance  \eqref{metr:vec:def}  of Section 3.  It reduces to the natural distance in the cone when $A$ and $B$ span the same subspace and it reduces to the Grassmann Riemannian distance between range$(A)$ and range$(B)$ when $A$ and $B$ are rank $p$ projectors.   Furthermore, it recovers a maximal number of the desirable invariance properties of the natural distance on the symmetric cone, exactly as expected from the planar example in Section 3. Homothetic transformations, isometries leaving the origin fixed, and pseudoinversion correspond to angle-preserving transformations (when $n>2$, these transformations are the conformal transformations of the Euclidean space). Note we already proved in Section 3 that the invariance with respect to transformations which do not preserve the angles is anyway an impossible property to obtain, at least when  $p$ is small enough.

The computation of the closeness \eqref{dist:def} involves the computation of principal angles and vectors, which is standard, and can be done via QR factorization at a numerical cost $O(np^2)$. The computation of the second term in \eqref{dist:def} involves a symmetric generalized eigenvalue problem, which requires $O(p^3)$ operations. The \emph{linear} complexity in the dimension $n$ makes the distance calculation efficient even in large-scale problems provided $p\ll n$.

The following proposition establishes a link between the closeness measure $l$ and the Riemannian distance $d_\Sy$ in $(\Sy,g)$.   
\begin{prop}\label{geodes:prop}Let $A=U_AR_A^2U_A^T$ and $B=U_BR_B^2U_B^T$ be two elements of $\Sy$. We have $$0\leq l^2(\gamma_{A \rightarrow B})-d^2_\Sy(A,B)\leq k~\max_{O\in O(p)}d^2_{\mathrm{P_p}}(R_B^2,O^TR_B^2O).$$\end{prop}The  right term is well defined for all $O\in O(p)$, and it is finite since $O(p)$ is a compact set. Moreover it tends to $0$ as $k\rightarrow 0$. Thus the curves defined above are always longer than the geodesics, but they yield good approximations of the geodesics when the Grassman contribution is highly penalized, i.e., $k$ is small (proof  in Appendix \ref{app2}).

We conclude this section by observing that horizontal geodesics provide a complete family of  line search curves in $\Sy$. The curve  emanating from $A=UR^2U^T$ in the direction
$(\Delta,D)\in\mathcal H_{(U,R^2)}$ admits the analytical characterization  \begin{align} 
\gamma_A(t)=U(t)R^2(t)U^T(t)\end{align}where $$U(t)=UV^T\cos{(\Gamma
t)}V+X\sin(\Gamma t)V$$ is the (Grassman geodesic) curve emanating
from $U$ in direction $\Delta$, i.e. $X\Gamma V=\Delta$ is the
compact SVD of $\Delta$; and $$R^2(t)=R\exp(tR^{-1}DR^{-1})R$$ is
the ($\mathrm{P_p}$ geodesic) curve emanating from $R^2$ in the
direction $D$. 
An interesting feature is that neither the curve $\gamma_{A}(\cdot)$ nor the geodesics reach the boundary of $\Sy$ in finite time (proof in Appendix \ref{app2}).
\begin{prop}\label{compl:prop}
For any $A$ in $\Sy$, the curves $\gamma_{A}(\cdot)$ can be extended from  $t=-\infty$ to $t=\infty$. Moreover, the manifold $(\Sy, g)$ is geodesically complete.
\end{prop}

\section{Geometric mean in \boldmath$\Sy$}\label{mean:sec}
{The   closeness measure between} $A$ and $B$ in
the previous section provides a direct formula for the {``halfway"} matrix
\begin{align}\label{mean:def}A\circ B=WKW^T,\end{align} where
$K=R_A(R_A^{-1}R_B^2R_A^{-1})^{1/2}R_A$ is the Riemannian mean  of
$R_A^2$ and $R_B^2$ in $\mathrm{P_p}$,
$W=\cos(\Gamma/2)U_A+\sin(\Gamma/2)X$ is the Riemannian mean of
range($A$) and range($B$), and \emph{where $R_A,R_B,U_A,X$ are
defined via the SVD \eqref{SVD}}.

The mean  is uniquely defined in the generic case  $\theta_{p}\neq\pi/2$, as the mean in Grassman between range($A$) and range($B$) is uniquely defined. When $\theta_{p}=\pi/2$ and $\theta_{p-1}\neq\pi/2$, according to Theorem 2 $l(\gamma_{A \rightarrow B})$ is uniquely defined, but there are two Riemannian means of range($A$) and range($B$) in Grassman. So the mean \eqref{mean:def} is well-defined but there are two midpoints. If $\theta_{p-1}=\pi/2$, there is an infinity of midpoints. Indeed, let $r$ be the number of principal angles equal to $\pi/2$. Looking at the proof of Theorem \ref{inv:dist:thm} in Appendix \ref{app2}, we see that the midpoint  will necessarily be a matrix of the form \eqref{mean:def} where $U_A,R_A$ are replaced by $U_AP,P^TR_AP$, with $P\in\mathrm{O(p)}$ representing an arbitrary rotation on the span of the $p-r$ last principal vectors. From now on we will systematically assume  that  $\theta_{p}\neq\pi/2$. 

The {geometric mean} $A\circ B$ is closely related to the popular geometric mean $A\#B$ of Ando \cite{ando78,pusz,ando}. One definition of Ando mean is to consider $A\#B$ as the solution of the extremal problem\begin{align}\label{extremal:eq}\max~ \left\{X\succeq0~:~\begin{pmatrix}A&X\\X&B\end{pmatrix}\succeq0 \right\}.\end{align}For full-rank  positive definite matrices $A$ and $B$, the solution of \eqref{extremal:eq} is given by  $X=A^{1/2}(A^{-1/2}BA^{-1/2})^{1/2}A^{1/2}$, that is, Ando mean coincides with the geometric mean \eqref{mean:def}. In contrast, for matrices $A$ and $B$ in $\Sy$, the geometric mean \eqref{mean:def} differs from the Ando mean \eqref{extremal:eq}.  In particular the geometric mean is by definition  \emph{rank-preserving}  whereas the solution of \eqref{extremal:eq} has a rank which is upper bounded by dim(range($A$) $\cap$ range($B$)).  A simple example is provided by the example already discussed in Section \ref{2D:sec}. The geometric mean of   $A$=diag(4,0) and
$B$=diag(0,1)  is $A\circ B=\begin{pmatrix}1&1\\1&1\end{pmatrix}$ while the solution of \eqref{extremal:eq} is $A\#B=\begin{pmatrix}0&0\\0&0\end{pmatrix}$. This is easy to check because the Ando mean can  be obtained by density: letting $A'=A+\epsilon I$ and $B'=B+\epsilon I$ the mean corresponds to the limit of $A'^{1/2}(A'^{-1/2}B'A'^{-1/2})^{1/2}A'^{1/2}$ for $\epsilon\rightarrow 0$.

Apart from this important difference, both the geometric mean \eqref{mean:def} and the Ando mean enjoy most desirable properties of a matrix geometric mean \cite{ando}, listed in the proposition below. 
\begin{prop}The mean $A\circ B$ possesses the  properties listed below.\end{prop}\begin{enumerate}
\item Joint homogeneity $\alpha A\circ \beta B=(\alpha\beta)^{1/2}A\circ B$.
\item Permutation invariance $A\circ  B=B\circ  A.$
\item Monotonicity. If $A\leq A_0$ (i.e. ($A_0-A$) is a positive matrix) and $B\leq B_0$, the means are comparable and verify $A\circ B\leq A_0\circ B_0$.
\item Congruence invariance. For any $(\mu,P)\in \RR\times \mathrm{O(n)}$ we have $(\mu P^TA\mu P)\circ (\mu P^TB\mu P)=\mu P^T(A\circ B )\mu P$.
\item Self-duality $(A\circ B)^\dag  =(A^\dag\circ B^\dag),$ where ``$\dag$" denotes pseudoinversion.
\end{enumerate}

The mean
\eqref{mean:def} may prove useful to generalize to low-rank positive semidefinite matrices the growing use of the geometric mean in applications requiring interpolation and the filtering of positive definite matrices; see e.g. \cite{pennec-06,arsigny,barbaresco,moakher06}.

A particular case of interest is the set of rank $p$ projectors \begin{align}\label{proj:eq}
\{P\in\RR^{n\times n}/~P^T=P,~P^2=P,~\tr{P}=p\},
\end{align}which is in bijection with the Grassman manifold of $p$-dimensional subspaces. Not surprisingly,  the geometric mean \eqref{mean:def} of $A$ and $B$ in \eqref{proj:eq} agrees with the Riemannian mean of range($A)$ and range$(B)$ on $\mathrm{Gr}(p,n)$.  On this set, the geometric mean also agrees with the Riemannian mean in $\Sy$ since $d_\Sy(A,B)= d_\mathrm{Gr(p,n)}(\text{range}(A),\text{range}(B))$ for $A$, $B$ projectors (see Appendix \ref{app2}).

On the other hand, when $A$ and $B$ have the same range, the geometric mean \eqref{mean:def} does not necessarily agree with the   Riemannian mean in $\Sy$ since the connecting curve $\gamma_{A \rightarrow B}(\cdot)$ is not necessarily  a geodesic. Note however that   the geometric mean is a good approximation of the Riemannian mean when $k$ is small. In fact the geometric mean \eqref{mean:def}  is a more  natural extension of the Ando mean than the Riemannian mean since   $A\circ B=U(R_A^2\#R_B^2)U^T$ for   $A=UR_A^2U^T$ and $B=UR_B^2U^T$. 

\section{Conclusion}\label{s7}
This paper generalizes the Riemannian geometry of the symmetric cone $\mathrm{P_n}$ to the manifold $\Sy$ of positive semidefinite matrices of fixed-rank $p$.
The generalization is based on the quotient geometry $\Sy=(\mathrm{V_{n,p}}\times \mathrm{P_p})/\mathrm{O(p)}$ that leads to a natural metric with decoupled
contribution in $\mathrm{Gr(p,n)}$ and the cone $\mathrm{P_p}$. This geometry leads to an explicit and natural   {notion of closeness in $\Sy$}, which in turn provides
a computable and natural definition of Riemannian distance and rank-preserving geometric mean.

The proposed computational tools may prove useful in applications 
involving computations with low-rank approximations of large-scale positive definite 
matrices. Such applications have already appeared in MRI tensor computing 
\cite{pennec-06,pennec2,fletcher,arsigny}, and in radar processing \cite{barbaresco}. 
Particular areas where the tools may prove useful in the future are the growing use 
of kernel-based methods, and low-rank approximations in machine learning and in 
bioinformatics. Kernel learning \cite{ratsch}, kernel completion \cite{scholkopf04}, 
and the use of Bregman divergence to address matrix nearness problems 
\cite{dhillon-07,kulis-06} (the natural distance on the cone is a Bregman divergence) 
are exemplary illustrations of areas that could benefit from the computational tools 
introduced in this paper.

\section*{Appendix}
\subsection{Reductive homogeneous space and invariant metrics\label{app2}} 
Let us choose the  identification $P_n\cong \mathrm{Gl(n)}/\mathrm{O(n)}$ as a 
starting point. Every matrix of the Lie algebra $\mathfrak{gl}(n)=T_I 
\mathrm{Gl(n)}=\RR^{n\times n}$ can be expressed as the sum of its symmetric part 
and its skew-symmetric part. Thus $\mathfrak{gl}(n)=\mathfrak m+\mathfrak{so}(n)$ 
is the direct sum of the symmetric matrices $\mathfrak m$ and the Lie algebra 
$\mathfrak{so}(n)=T_I \mathrm{O(n)}$, made of skew-symmetric matrices. Moreover, 
for any $S\in \mathfrak m$ and $O\in \mathrm{O(n)}$ we have 
$Ad_O(S)=OSO^{-1}=OSO^T\in\mathfrak m$. Thus $Ad_{\mathrm{O(n)}}(\mathfrak m)\subset \mathfrak m$. 
The existence of $\mathfrak m$ such that $\mathfrak{gl}(n)=\mathfrak m+\mathfrak{so}(n)$ 
is a direct sum, and $Ad_{\mathrm{O(n)}}(\mathfrak m)\subset \mathfrak m$ proves that 
$\mathrm{Gl(n)}/\mathrm{O(n)}$ can be called a \emph{reductive} homogeneous space.

We provide some more information for the readers who are familiar
with standard ideas of Lie group theory. Let $\{Z(t),~t>0\}$ be a
trajectory in $\mathrm{Gl(n)}$. $\{Z_1(t)=Z(t)O(t),~t>0\}$ with
$O(t)\in \mathrm{O(n)}$ represents the same trajectory in $P_n$.
Let $\omega(t)=Z^{-1}\dotex Z\in\mathfrak{gl}(n)$ and
$\omega_1(t)=Z_1^{-1}\dotex Z_1$. We have
$\omega(t)=Ad_{O(t)}\bigl(\omega_1(t)-O(t)^{-1}\dotex O(t)\bigr)$
(see~\cite{klein}). Decomposing $\omega,\omega_1$ in their
symmetric and skew-symmetric parts, we see that any
$Ad_{\mathrm{O(n)}}$-invariant scalar product on $\mathfrak m$
will allow the construction of a well-defined
$\mathrm{GL(n)}$-invariant metric on the homogeneous space
$\mathrm{Gl(n)}/\mathrm{O(n)}$ corresponding to this scalar
product at $I$. All  $Ad_{\mathrm{O(n)}}$-invariant scalar
products on the symmetric matrices are given up to a constant
scale factor by $\tr{X_1^TX_2}+\beta\tr{X_1}\tr{X_2}$ with
$\beta>-\frac{1}{n}$. Indeed, they are derived from rotationally
invariant norms $\norm{X}^2$ on symmetric matrices. Thus they can
only depend on the scalar  invariants $\tr{X}$, $\tr{X^2}$,
$\tr{X^3}$, etc. As they are quadratic functions they can only
depend on $\tr{X}^2$ and $\tr{X^2}$. The condition on $\beta$
ensures positive definiteness (see, e.g.,~\cite{pennecHDR}).

\subsection{Several proofs\label{app2}}\mbox{}\paragraph{Proof of Theorem \ref{inv:dist:thm}}

Let $U(t)$ be a geodesic linking $U_A$ and $U_B$ in $\mathrm{V_{n,p}}$  and $R^2(t)$ a geodesic linking $R_A^2$ and $R_B^2$ in $\mathrm{P_p}$, where $(U_A,R_A^2)$ and $(U_B,R_B^2)$  are representatives of $A$ and $B$ defined via SVD \eqref{SVD}. By corollary 3.57 of \cite{Oneill-book} $(U(t),R^2(t))$ is a geodesic curve in the structure space $\mathrm{V_{n,p}}\times \mathrm{P_p}$ endowed with the metric $g^{SP}$. The choice of the representatives is such that for every $t>0$,   $\dotex U(t)$ is orthogonal to $U(t)$, as it is an element of the horizontal space of  $\mathrm{Gr(p,n)}=\mathrm{V_{n,p}}/\mathrm{O(p)}$. This has two consequences. First the length of the geodesic in the structure space is 
$$d^2(A,B)={d^2_{\mathbb Gr(p,n)}(\text{range}(A),\text{range}(B))+k~d^2_{\mathrm{P_p}}(R_A^2,R_B^2) }.$$
Because  the Riemannian distance in $\mathrm{Gr(p,n)}$ is $\norm{\Theta}_F$ \cite{EAS98,AMS2004-01}  and  the Riemannian distance in $\mathrm{P_p}$ is $\norm{\log R_A^{-1}R_B^2 R_A^{-1}}_F$, we obtain 
$$d^2(A,B)=\norm{\Theta}_F^2+k~\norm{\log R_A^{-1}R_B^2 R_A^{-1}}_F^2$$
Then $\dotex(U(t),R^2(t))$ belongs to the horizontal space $\mathcal H_{(U(t),R^2(t))}$ for all $t>0$, that is, $(U(t),R^2(t))$ is a {\it horizontal} curve. As $g^{SP}$ induces the metric \eqref{metric:def} on the quotient manifold, the length $l(\gamma_{A \rightarrow B})$ of  $\gamma_{A \rightarrow B}$ in the quotient space $(\mathrm{V_{n,p}}\times \mathrm{P_p})/\mathrm{O(p)}$ is also $d(A,B)$, proving 
\eqref{dist:def}.

{Uniqueness.}
First of all suppose the $\sigma_i$'s are distinct. According to the uniqueness of the SVD in \eqref{SVD}, the matrices $U_A$ and $U_B$ are unique (up to a joint multiplication of any columns by $-1$, which is an orthogonal transformation) and do not depend on the choice of $V_A$ and $V_B$. Under the less restrictive assumption $\sigma_{p-1}>0$, suppose there exists $1\leq i \leq p-1$ such that $\sigma_i=\sigma_{i+1}$. Then $\sigma_i>0$ and  the SVD yields $u_i^A,u_i^B,u_{i+1}^A,u_{i+1}^B$ such that $(u_i^A)^Tu_i^B=(u_{i+1}^A)^Tu_{i+1}^B=\sigma_i$. For the sake of simplicity we assume there are only $2$ principal vectors associated to $\sigma_i$. The generalization to an arbitrary number of vectors is straightforward and leads to the same conclusion. The SVD yields non-unique principal vectors, since any normalized linear combination of $u_{i}$ and $u_{i+1}$ is still a principal vector. Let $u^A=au_i^A+bu_{i+1}^A$ and $u^B=cu_i^B+du_{i+1}^B$ be other principal vectors associated to $\sigma_i$. By definition $(u^A)^Tu^B=\sigma_i$ and $(u^A)^Tu^A=(u^B)^Tu^B=1$. Since $(u_i^A)^Tu_{i+1}^B=0$  we have $(ac+bd)\sigma_i=\sigma_i$. It implies $a=c$ and $b=d$ (Cauchy-Schwartz equality) and, necessarily, $u^A$ and $u^B$ are obtained from $u^A_i,u^A_{i+1}$, and $u^B_i,u^B_{i+1}$ via the same orthogonal transformation (unless $\sigma_i=0$, which is impossible). We proved the principal vectors $U_A$ and $U_B$ are defined up to a \emph{joint} multiplication on the right by a block orthogonal matrix $P\in \mathrm{O(p)}$, each block being associated to the same eigenvalue. Let $U_AP=U_A'$ and $U_BP=U_B'$, and let $A=U_A'R_A'^2(U_A')^T$ and $B=U_B'R_B'^2(U_B)'^T$. We have $
R_A'^2=P^TR_A^2P,~R_B'^2=P^TR_B^2P$ and $l(\gamma_{A \rightarrow B})$ is unchanged
since $d_{\mathrm{P_p}}(R_A^2,R_B^2)=d_{\mathrm{P_p}}(R_A'^2,R_B'^2)$. Note that the reciprocal of
this result is also true.

{Invariances}. Using the preceding paragraph, we know that  \eqref{dist:def} is well defined under the basic assumption $\sigma_{p-1}>0\Leftrightarrow \theta_{p-1}<\pi/2$. Thus the following features are sufficient to complete the proof:\begin{itemize}
\item For $\mu\in\RR$ we have $\mu^2 A=U_A(\mu^2 R_A^2)U_A^T$. Since  we have $d_{\mathrm{P_p}}(\mu^2R_A^2,\mu^2R_B^2)=d_{\mathrm{P_p}}(R^2_A,R^2_B)$  we have $l(\gamma_{\mu^2A \rightarrow \mu^2B})=l(\gamma_{A \rightarrow B})$.
\item Let $O\in \mathrm{O(n)}$, and let $U_A$ and $U_B$ be the principal vectors associated to $A$, $B$. Then $U_A'=OU_A$ and $U_B'=OU_B$ are principal vectors associated to $OAO^T$ and $OBO^T$ since orthogonal transformations preserve the angles: $(OU_A)^TOU_B=U_A^TU_B$, so the Grassman distance is unchanged. Since $OAO^T=U_A'R^2_AU'^T~\text{and}~OBO^T=U_B'R^2_BU_B'^T$ we see that $R^2_A,R^2_B$  are also unchanged by the transformation and $l(\gamma_{OAO^T \rightarrow OBO^T})=l(\gamma_{A \rightarrow B})$.
\item$A^\dag=U_AR^{-2}_AU_A^T$ is the pseudoinverse of $A$.  $l(\gamma_{A \rightarrow B})-l(\gamma_{A^\dag \rightarrow B^\dag})=d_{\mathrm{P_p}}(R^{-2}_A,R^{-2}_B)-d_{\mathrm{P_p}}(R_A^2,R_B^2)=0.$\end{itemize}

\paragraph{Proof of Proposition 2} 
First, as $\gamma_{A\rightarrow B}$ induces a non-minimal path connecting $A$ and $B$ in $\Sy$ we have $l^2(\gamma_{A \rightarrow B})\geq d_\Sy(A,B)$. Then, let $\gamma(t)$ be a geodesic linking $A$ and $B$ in $\Sy$, with $\gamma(0)=A, \; \gamma(1)=B$ and $d_{\Sy}(A,B)$ its length. As proved in the sequel $\gamma(t)$ can be globally horizontally lifted in the structure space $\mathrm{V_{n,p}}\times \mathrm{P_p}$. The horizontal lift $\tilde\gamma(t)$ at $(U_A,R_A)$  is a horizontal curve in the structure space  $\mathrm{V_{n,p}}\times \mathrm{P_p}$ connecting the two fibers. Therefore there exists $O\in O(p)$ such that $\tilde\gamma(1)= (U_BO,O^TR_B^2O)$, where $(U_B,R^2_B)$ is defined via the SVD \eqref{SVD}. The length of $\tilde\gamma$ in the structure space is also $d_{\Sy}(A,B)$. Since $\tilde\gamma$ is not necessarily a minimal curve in the structure space we have  $d^2_{\Sy}(A,B)\geq d^2_{\mathrm{V_{n,p}}}(U_A,U_BO)+k ~d^2_{\mathrm{P_p}}(R_A^2,O^TR_B^2O)$. 
Because  distances in  $\mathrm{ Gr}(p,n)$ are shorter than those in $\mathrm{V_{n,p}}$, this implies
$d^2_{\Sy}(A,B) \geq d^2_{\mathrm{ Gr}(p,n)}(\text{range}(U_A),\text{range}(U_BO))+k ~d^2_{\mathrm{P_p}}(R_A^2,O^TR_B^2O)$.   Finally we proved 
\begin{eqnarray*}
l^2(\gamma_{A \rightarrow B}) & \leq  & d^2_{\mathrm{ Gr}(p,n)}(\text{range}(U_A),\text{range}(U_B))+k ~d^2_{\mathrm{P_p}}(R_A^2,O^TR_B^2O)+k ~d^2_{\mathrm{P_p}}(O^TR_B^2O,R_B^2) \\
& \leq &  d^2_{\Sy}(A,B)+k ~d^2_{\mathrm{P_p}}(R_B^2,O^TR_B^2O)
\end{eqnarray*}

\paragraph{Proof of Proposition 3}The structure space is geodesically complete (it is a product of two geodesically complete Riemannian manifolds whose  geodesics are product of geodesics). Thus horizontal geodesics are defined for all $t\in\RR$.

A manifold is geodesically complete if and only if it is complete as a metric space (see Theorem 10.3 Ch I in \cite{Hegalson}). Thus the structure space is complete. To prove that the quotient space $\Sy$ is geodesically complete it suffices to prove that it is complete. Let $\{u_m\}$ be a Cauchy sequence in the quotient space $\Sy$. One can find a broken geodesic $\gamma(t)$ linking the $u_i$'s such that $\gamma(t_i)=u_i$. Consider an horizontal lift $\tilde\gamma$ of $\gamma$ (it is well defined as proved below). The distance between  $\tilde\gamma(t_i)$ and $\tilde\gamma(t_j)$ in the structure space is  $d_\Sy(\gamma(t_i),\gamma(t_j))= d_\Sy(u_i,u_j)$  and thus $\{\tilde\gamma(t_m)\}$ is a cauchy sequence in the structure space. Thus it converges to a point $\tilde p=(U,R^2)$. Let us prove $\{u_m\}$ converges to $p=UR^2U^T$. Thanks to the exponential map in the structure space, $\tilde\gamma(t_i)$  can be identified to an element $\xi_i=(\Delta_i+U\Omega_i,D_i)\in T_{(U,R^2)}\mathrm{V_{n,p}}\times \mathrm{P_p}$ for $i$ large enough. The geodesic distance in the structure space between $\tilde\gamma(t_i)$ and $\tilde p$ is $g^{SP}(\xi_i,\xi_i)=\norm{\Omega_i}^2+\norm{\Delta_i}^2+k~\tr{(D_iR^{-2})^2}$ and thus $\norm{\Omega_i}\rightarrow 0$. Thus the coordinate of $\tilde\gamma(t_i)$ along the fiber tends to zero so $d_\Sy(U_iR_i^2U_i^T,UR^2U^T)\rightarrow 0$.

Any geodesic $\gamma$ of $\Sy$ can be globally horizontally lifted. This result was shown by R. Hermann when the horizontal and vertical spaces are orthogonal. The result holds in our case and the proof is unchanged. Let us recall the main points of the proof : let $\tilde\gamma_p$ is the maximal horizontal lift of $\gamma$ at $p$, defined on, say $[0,\epsilon[$. It is enough to show that $\tilde\gamma_p$ can be extended at $\epsilon$. The continuity of $\gamma$  and the completeness of the structure space imply the existence of a limit point $q=\lim_{t\rightarrow\epsilon}\tilde\gamma_p(t)$.

\subsection*{Acknowledgment}
The authors wish to thank P.-A. Absil for several discussions on the paper.


\end{document}